\documentclass{amsart}
\usepackage{amssymb, amsmath}
\usepackage{bbm}
\usepackage{tikz-cd}

\title[MES Preservers]{On the Preservers of Maximally Entangled States}
 \author[B. W. Grossmann and H. J. Woerdeman]{Ben Grossmann}
\author[]{Hugo J.~Woerdeman}
\address{Department of Mathematics \\
Drexel University\\
3141 Chestnut St.\\
Philadelphia, PA, 19104}
\email{benwgrossmann@gmail.com,hugo@math.drexel.edu}

\thanks{HJW is supported by Simons Foundation grant 355645.}


\date{}

\newtheorem{thm}{Theorem}[section]

\newtheorem{rmk}[thm]{Remark}

\newtheorem{lem}[thm]{Lemma}

\newtheorem{prop}[thm]{Proposition}

\numberwithin{equation}{section}

\newcommand{\proofend}{\hfill $\Box$ \\[.25in] }

\newcommand{\C}{\mathbb{C}}
\newcommand{\R}{\mathbb{R}}

\newcommand{\Q}{\mathbb{Q}}

\newcommand{\s}[1]{\mathcal{#1}} 

\newcommand{\mb}[1]{\mathbf{#1}}

\DeclareMathOperator{\ad}{ad}
\DeclareMathOperator{\tr}{tr}
\DeclareMathOperator{\vecop}{vec}

\DeclareMathOperator{\MES}{MES}
\DeclareMathOperator{\spn}{span}

\DeclareMathOperator{\coisom}{coisom}
\DeclareMathOperator{\diag}{diag}

\DeclareMathOperator{\re}{Re}

\newcommand{\id}{\mathbbm{1}}

\newcommand{\mat}[1]{
\begin{bmatrix} #1
\end{bmatrix}}

\newcommand{\smat}[1]{
\left[\begin{smallmatrix} #1
\end{smallmatrix}\right]}

\allowdisplaybreaks

\begin{document}

\maketitle

\begin{abstract}
We characterize the linear maps that preserve maximally entangled states in $L(\s X \otimes \s Y)$ in the case where $\dim(X)$ divides $\dim(Y)$.
\end{abstract}

\section{Introduction}


Entanglement is considered a valuable resource in quantum information theory; it is responsible for the power of quantum computing and for applications such as superdense coding and
quantum teleportation.  
In the interest of ``conserving'' the resource of entanglement, we might ask the following question: which linear maps produce a maximally entangled output for any maximally entangled input?

More concretely: a quantum system with $m$ mutually exclusive configurations can be represented by an $m$-dimensional inner product space, which we will call $\s X$.  The configuration of such a system is represented with a unit vector $u \in \s X$.  If $\s X, \s Y$ are the spaces corresponding to two such systems, then a configuration of the joint system can be represented with a unit vector $u \in \s X \otimes \s Y$, where $\s X \otimes \s Y$ is the tensor product of the Hilbert spaces associated with each system.  We could also represent this system with the associated ``pure state'' density operator $A = uu^*$.  We can measure the ``entropy of entanglement'' of $A$ by $S(\tr_{\s Y}(A))$, where $\tr_{\s Y}$ denotes a partial trace over the space $\s Y$ and $S(\rho) = -\tr(\rho \log_2 \rho)$ is the Von Neumann entropy.  We refer to the set of pure state density operators that maximize this measure of entanglement ``maximally entangled states,'' and denote the set as $\MES_{\s X, \s Y}$.  This set can be succinctly described as
\begin{equation*}
    \MES_{\s X, \s Y} = \left\{\frac 1m \vecop(A)\vecop(A)^* \mid A : \s Y \to \s X \text{ is a coisometry}\right\},
\end{equation*}
where $\vecop$ stacks the rows of a matrix $A$ (see Section 2 for details). Quantum information is generally transmitted via a ``quantum channel''.  Mathematically, this channel is a linear map that takes density matrices over a first system to density matrices over a second system (with an additional ``complete positivity'' condition).  With this in mind, we hope to classify the linear maps $\Phi:\spn(\MES_{\s X, \s Y}) \to \spn(\MES_{\s X, \s Y})$ that preserve the property of maximal entanglement, which is to say that they satisfy $\Phi(\MES_{\s X, \s Y}) \subset \MES_{\s X, \s Y}$.

The problem of classifying the linear maps that ``preserve'' the subset $\MES_{\s X, \s Y}$ is an example of a ``linear preserver problem''. Linear preserver problems like this one have a long history (see, for instance, \cite{PresSurvey2000} and \cite[Chapter 7]{FuncId}).  Recently, there has been work done on linear preservers within the context of quantum information theory. Two relevant papers to this endeavor are \cite{QIAut} and the more recent paper \cite{Poon15}, whose results we generalize here.


In this paper, we prove the following result.

\begin{thm}\label{pres}
Let $\s Y = \s X^k := \s X \oplus \cdots \oplus \s X$, $k>1$. If $\Phi:\spn(\MES_{\s X,\s Y}) \to \spn(\MES_{\s X, \s Y})$ is an invertible linear map that preserves $\MES_{\s X,\s Y}$, then $\Phi$ is of the form
\[
\Phi: X \mapsto (U \otimes V)\,X^\sigma\,(U \otimes V)^*
\]
where $U\in L(\s X),V \in L(\s Y)$ are unitary if $A \mapsto A^\sigma$ denotes either the identity or transpose map.
\end{thm}


This generalizes (the invertible case of) Poon's result from \cite{Poon15}, which can be stated within our framework as follows:

\begin{thm} \label{Poonthm1}

\cite{Poon15} A linear map $\Phi: \spn(\MES_{\s X,\s X}) \to \spn(\MES_{\s X, \s X})$ preserves $\MES_{\s X, \s X}$ if and only if it has one of the following forms:
\begin{enumerate}
    \item $\Phi(A \otimes B) = (U \otimes V)(A \otimes B)^\sigma(U \otimes V)^*$
    \item $\Phi(A \otimes B) = (U \otimes V)(B \otimes A)^\sigma (U \otimes V)^*$
    \item $\Phi(X) = (\tr X) \rho$ for some $\rho \in \MES$
\end{enumerate}
where $A \mapsto A^\sigma$ denotes either the identity or transpose map.
\end{thm}

Notice that (2) does not appear in Theorem \ref{pres}.  Indeed, for (2) to make sense one needs $\s X = \s Y$.  In addition, Poon was able to describe the non-invertible preservers, which take the form (3).


The proof can be summarized as follows.  First, we show that an MES-preserving map $\Phi$ induces a continuous map $\zeta$ over the projective space $P(\coisom(\s Y, \s X))$ (to be introduced in Section 4).  Moreover, we show that for mutually orthogonal coisometries $A,B$, there exist mutually orthogonal coisometries $C,D$ such that we have either
\begin{align*}
    \zeta([\alpha A + \beta B]) = [\alpha C + \beta D] \quad \text{or} \quad \zeta([\alpha A + \beta B]) = [\bar \alpha C + \bar \beta D].
\end{align*}
By considering the degree of $\zeta$ restricted to the span of $[A]$ and $[B]$, we deduce whether the $A \mapsto A^\sigma$ of our $\MES$-preserving $\Phi$ is the identity or the transpose. 

Having classified the possible maps $\zeta$, we then construct an $\MES$-preserving extension $\tilde \Phi : \MES_{\s Y, \s Y}  \to \MES_{\s Y, \s Y}$ and apply Theorem \ref{Poonthm1} to this extension $\tilde \Phi$. By considering the restriction of this map $\tilde \Phi$ to a subspace identified with $\MES_{\s X, \s Y}$, we conclude that the restriction $\Phi$ must be of the form described in Theorem \ref{pres}.


The paper is organized as follows.  In Section 2, we introduce definitions, particularly the definition of $\MES$ (the set of maximally entangled states) and the definition of our vectorization operator.  In Section 3, we discuss some characterizations of $\MES_{\s X, \s Y}$, and how it relates to the set of coisometries $\coisom(\s Y, \s X)$.  In Section 4, we define the continuous maps on $\MES$ that we will rely on in order to define the extension $\tilde \Phi$.  In Section 5, we consider how $\zeta$ behaves on pairwise-orthogonal coisometries.  In Section 6, we use $\zeta$ to define the extension $\tilde \Phi:\MES_{\s Y, \s Y} \to \MES_{\s Y, \s Y}$.  Finally, in Section 7 we prove Theorem \ref{pres}, our main result.

\section{Definitions}

Let $\s X, \s Y$ denote (finite dimensional, complex) Hilbert spaces with $\dim(\s X) = m$, $\dim(\s Y) = n$.  Let $\s X \otimes \s Y$ denote the tensor product of these spaces, which is itself a Hilbert space of dimension $mn$.  Let $L(\s X, \s Y)$ denote the space of linear transformations between $\s X$ and $\s Y$.  We will also use the shorthand $L(\s X) = L(\s X, \s X)$.  We will use $\id_{\s X}$ to denote the identity operator over the Hilbert space $\s X$.

We say that an operator $A: \s X \to \s Y$ is an \textit{isometry} if it satisfies $A^*A = \id_{\s X}$.  We will say that $A$ is a \textit{coisometry} if its adjoint operator $A^*: \s Y \to \s X$ is an isometry.

A unit vector $u \in \s X \otimes \s Y$ is said to be \textit{maximally entangled} if it can be written in the form
\begin{equation}
u = \sum_{j=1}^r \frac{1}{\sqrt{r}} x_j \otimes y_j
\end{equation}
for orthonormal sets of vectors $\{x_1,\dots,x_r\} \subset \s X,\{y_1,\dots,y_r\} \subset \s Y$ and $r = \min\{m,n\}$.

The \textit{pure state} (rank-$1$, trace-$1$, positive semidefinite operator) associated with $u$ is the projection operator $uu^* \in L(\s X \otimes \s Y)$.  Such an operator is called a \textit{maximally entangled state}.  That is, we define the set $\MES_{\s X,\s Y}$ of maximally entangled states over $\s X$ and $\s Y$ to be
\begin{equation}
    \MES_{\s X, \s Y} = \{uu^* : u \text{ is maximally entangled}\} \subset L(\s X \otimes \s Y).
\end{equation}
We will sometimes denote $\MES_{\s X, \s Y}$ as $\MES$ (without its subscripts) when the context is clear.

We define the vectorization operator $\vecop:L(\s Y, \s X) \to \s X \otimes \s Y$ to be the linear map satisfying
\begin{equation}
    \vecop(xy^T) = x \otimes y \qquad \text{for all } x \in \s X,\ y \in \s Y.
\end{equation}
Notably, $\vecop$ is an isomorphism of Hilbert spaces.  It may be helpful to note that, in terms of matrices, $\vecop$ can be written as $\vecop:\C^{m \times n} \to \C^{mn}$, 
\begin{equation}
    \vecop({[a_{ij}]_{i=1}^m}_{j=1}^n) = (a_{11}, \dots, a_{1n},a_{21},\dots,a_{2n},\dots,a_{m1}, \dots, a_{mn})^T ,
\end{equation}
which is to say that our vectorization operator ``stacks the rows'' of a matrix (in contrast to the other common convention of ``stacking the columns'').  We adopted this convention from \cite{TQI}.

Throughout the paper $e_1,\dots,e_n$ will denote the canonical basis for $\C^n$ (for instance, $e_1 = (1,0,\dots,0)^T$).  We will also take $E_{ij} \in \C^{m \times n}$ to be the matrix $E_{ij} = e_ie_j^T$, i.e. the matrix that has a $1$ as its $i,j$ entry and zeros elsewhere.

\begin{rmk}
In the remainder of the paper, we will take $\dim(\s X) = m < \infty$ and $\dim(\s Y) = n < \infty$, with $m \leq n$.
\end{rmk}

\section{Characterizing $\MES$}

\begin{prop} \label{prop1}
If $n \geq m$, $u$ is maximally entangled in $\s X \otimes \s Y$ iff there exists a coisometry $A:\s Y \to \s X$ such that $u = \frac 1{\sqrt{m}} \vecop(A)$.

\end{prop}

\textit{Proof:} Suppose that $u$ is a maximally entangled unit vector.  Let $\{x_1,\dots,x_m\}$ and $\{y_1,\dots,y_n\}$ be orthonormal bases of $\s X$ and $\s Y$ such that\\ 
$u = \frac {1}{\sqrt{m}}\sum_{j=1}^m x_j \otimes y_j$.  It follows that
\begin{equation}
    \vecop^{-1}(u) = \frac 1{\sqrt m} \sum_{j=1}^m \vecop^{-1}(x_j \otimes y_j) = \frac 1{\sqrt{m}} \left(\sum_{j=1}^m x_j y_j^T\right)
\end{equation}
and we recognize that $A = \sum_{j=1}^m x_j y_j^T$ is a singular value decomposition of a coisometry.

Conversely, suppose $u = \frac 1{\sqrt m} \vecop(A)$.  Then $A$ has a singular value decomposition $A = \sum_{j=1}^m x_j y_j^T$
so that 
\begin{equation}
    u = \frac 1{\sqrt m} \vecop \left(\sum_{j=1}^m x_j y_j^T\right) 
    = \frac 1{\sqrt m} \sum_{j=1}^m x_j \otimes y_j
\end{equation}
is maximally entangled. \proofend

Let $\tr_{\s Y}: L(\s X \otimes \s Y) \to L(\s X)$ denote the ``partial trace'' over the space $\s Y$.  That is, we define $\tr_{\s Y}$ to be the unique linear operator satisfying
\begin{equation}
    \tr_{\s Y}(A \otimes B) = \tr(B) A \qquad \text{for all }A \in L(\s X), \quad B \in L(\s Y).
\end{equation}

The following lemma is stated, for instance, in Equation 1.133 of \cite{TQI}.

\begin{lem} \label{PT}
For any operators $A,B:\s Y \to \s X$, we have
\begin{equation}
    \tr_{\s Y}(\vecop(A)\vecop(B)^*) = AB^*.
\end{equation} 

\end{lem}

\textit{Proof:} Suppose first that $A$ and $B$ are rank-$1$.  We can then write $A = uv^T$ and $B = xy^T$ for $u,x \in \s X$ and $v,y \in \s Y$.  We find that
\begin{align*}
    \tr_{\s Y}(\vecop(A)\vecop(B)^*) &= 
    \tr_{\s Y}((u\otimes v)(x \otimes y)^*) 
    \\ & = 
    \tr_{\s Y}((ux^*) \otimes (vy^*))
    \\ & = (y^*v) ux^* = u\overline{v^*y}x^* = (uv^T)(xy^T)^* = AB^*.
\end{align*}
Now, if $A = \sum_j A_j$ and $B = \sum_k B_k$ where $A_j,B_k \in L(\s Y, \s X)$ have rank $1$, then 
\begin{align*}
    \tr_{\s Y}(\vecop(A)\vecop(B)^*) &= 
    \tr_{\s Y}\left(\vecop\left(\sum_{j}A_j\right)\vecop\left(\sum_k B_k\right)^*\right)
    \\ & = 
    \sum_{j,k} \tr_{\s Y} (\vecop(A_j)\vecop(B_k)^*)
    = \sum_{j,k} A_jB_k^* 
    \\ & = 
    \left(\sum_{j}A_j\right)\left(\sum_k B_k\right)^* = AB^*.
\end{align*}
Since any operator can be written as a sum of rank-1 operators, the desired conclusion holds. 

\proofend

\begin{lem} \label{PTMES}

For all $M \in \MES_{\s X, \s Y}$, we have $\tr_{\s Y}(M) = \frac 1m \id_{\s X}$.  Moreover, for $M \in \spn(\MES_{\s X, \s Y})$, we have $\tr_{\s Y}(M) = \frac{\tr(M)}{m}\id_{\s X}$.

\end{lem}

\textit{Proof:} By Proposition \ref{prop1}, there exists a coisometry $A: \s Y \to \s X$ such that $M = \frac 1m \vecop(A)\vecop(A)^*$.  By Lemma \ref{PT}, we find that
\begin{equation}
    \tr_{\s Y}(M) = \frac 1m \tr_{\s Y}(\vecop(A)\vecop(A)^*) = 
\frac 1m AA^* = \frac 1m \id_{\s X}.
\end{equation}
By the linearity of the partial trace, it follows that $\tr_{\s Y}(M) = \frac{\tr(M)}{m}\id_{\s X}$ for all $M \in \spn(\MES_{\s X,\s Y})$.

\proofend

The following proposition generalizes Remark 2.3 in \cite{Poon15}.

\begin{prop} \label{MESspan}

The set of pure states in the real linear span of MES is precisely MES.  That is: if $M$ is a rank-1, trace-1, positive operator, then $M \in \spn(\MES)$ implies that $M \in \MES$.

\end{prop}

\textit{Proof:} 
Let $M$ be a rank-1, trace-1, positive operator in $\spn(\MES)$.  Let $u$ be a unit-vector such that $M = uu^*$.  Let $A: \s Y \to \s X$ be an operator such that $u = \frac 1{\sqrt m}\vecop(A)$.  By Lemma \ref{PTMES}, $M \in \spn(\MES_{\s X, \s Y})$ implies that
\begin{equation}
    \tr_{\s Y}(M) = \frac {\tr M}{m}\id_{\s X} = \frac 1m \id_{\s X}.
\end{equation}
By Lemma \ref{PT}, this implies that
\begin{equation}
    \id_{\s X} = m \tr_{\s Y}(M) = \tr_{\s Y}(\vecop(A)\vecop(A^*)) = AA^*.
\end{equation}

That is, $A$ is a coisometry. By Proposition \ref{prop1}, we may conclude that $u = \frac 1{\sqrt m}\vecop(A)$ is maximally entangled, and that $M = uu^*$ is an element of $\MES$. 

\proofend

\section{Continuous Maps on $\MES$}

\begin{rmk}
In the remainder of the paper, $\Phi$ will denote an invertible linear map $\Phi: \spn(\MES_{\s X,\s Y}) \to \spn(\MES_{\s X,\s Y})$ such that $\Phi(\MES_{\s X, \s Y}) \subset \MES_{\s X, \s Y}$.
\end{rmk}

For a (finite-dimensional) vector space $V$ and a subset $S \subset V$, let $P(S)$ denote complex projective space over $S$, the topological space of one-dimensional $\C$-subspaces of $V$ that contain a non-zero element of $S$. For an element $A \in S$, we use $[A]=\{ z A : z \in {\mathbb C} \}$ to denote the element of $P(S)$ containing $A$.

For $A \in L(\s Y, \s X)$, let $\pi_A \in L(\s X \otimes Y)$ be given by
\[
\pi_A = \frac 1{\tr(AA^*)} \vecop(A)\vecop(A)^*.
\]
We see that $\pi_A \in \MES_{\s X, \s Y}$ if and only if $A: \s Y \to \s X$ is a (complex) multiple of a coisometry (that is, if $AA^* = k\id_{\s X}$ for some $k > 0$).

\begin{rmk}
Let $K$ be a compact subset of the finite-dimensional normed vector space $L(\s Y, \s X)$. The map $\pi$ is constant over any equivalence class in $P(K)$, so the induced map from $P(K)$ to $L(\s X \otimes \s Y)$ given by
$[A] \mapsto \pi_{A}$
is well-defined. Moreover, this induced map is a homeomorphism between $P(K)$ and the map's image.  
\end{rmk}

With $K = \coisom(\s Y, \s X)$ in particular, we see that the image of the map $[A] \mapsto \pi_A$ is precisely $\MES_{\s X, \s Y}$. 

\vspace{10 pt}

\textit{Proof:}
To see that the map $[A] \mapsto \pi_{A}$ is well defined and continuous over the quotient space $P(K)$, it suffices to note that if $A_1,A_2 \in K$ satisfy $A_1 \sim A_2$ (i.e. if $A_1 = k A_2$ for some non-zero $k\in \C$), then $\pi_{A_1} = \pi_{A_2}$.  We see moreover that this induced map is injective: that is, we observe that if $\pi(A_1) = \pi(A_2)$, then $A_1 \sim A_2$.  Because $K$ is compact, the quotient space $P(K) = K/\sim$ is also compact.

So, we have found that the induced map $[A]\mapsto \pi_A$ is a continuous, injective map from the compact quotient space $K/\sim$ to the Hausdorff space $L(\s X \otimes \s Y)$.  Because our map is an injective, continuous map with compact domain and Hausdorff codomain, it is a homeomorphism between the map's domain and the map's image (by Theorem 26.6 of \cite{Munkres}), as desired.

\proofend

With the above in mind, we note that the set $\coisom(\s Y, \s X)$ is compact in $L(\s Y, \s X)$.  We define $\zeta:P(\coisom(\s Y, \s X)) \to P(\coisom(\s Y, \s X))$ to be the unique linear map such that for all $A \in \coisom(\s Y, \s X)$, we have $\pi_{\zeta[A]} = \Phi(\pi_A)$.

That is, we define $\zeta$ so that the following diagram (of continuous maps) commutes:

\begin{equation}
\begin{tikzcd}[row sep=huge]
\MES_{\s X, \s Y} \arrow[r,"\Phi"]  
&\MES_{\s X,\s Y}\\
P(\coisom(\s Y, \s X)) \arrow[u,"\pi"] 
\arrow[r,dashed,"\zeta"]
& P(\coisom(\s Y, \s X)) \arrow[u,"\pi"]
\end{tikzcd}
\end{equation}

Since $\Phi$ is a continuous, invertible map on $\MES_{\s X,\s Y}$ and $\pi$ is a homeomorphism, we see that $\zeta$ must also be a continuous, invertible map on $P(\coisom(\s Y, \s X))$.
%
%


This map $\zeta$ will be particularly important to the construction of $\tilde \Phi$ in Section 6.  To begin, we show that $\zeta$ satisfies the following ``preservation of subspaces'' property.

\begin{prop} \label{3}
For $A_1,A_2,A_3 \in \coisom(\s Y, \s X)$: if $[A_3]\subset [A_1] + [A_2]$, then  $\zeta[A_3]$ must be a subspace of $\zeta[A_1]+\zeta[A_2]$.
\end{prop}

\textit{Proof:} Equivalently, we wish to show that
\begin{align*}
\exists \alpha > 0 &\text{ such that } \pi_{A_1} + \pi_{A_2} - \alpha \pi_{A_3} \text{ is positive semidefinite} \implies\\
\exists \beta > 0 &\text{ such that } \pi_{\zeta[A_1]} + \pi_{\zeta[A_2]} - \beta \pi_{\zeta[A_3]} \text{ is positive semidefinite}.
\end{align*}
That is, we wish to show that
\begin{align*}
\exists \alpha > 0 &\text{ such that } \pi_{A_1} + \pi_{A_2} - \alpha \pi_{A_3} \text{ is positive semidefinite} \implies\\
\exists \beta > 0 &\text{ such that } \Phi(\pi_{A_1} + \pi_{A_2} - \beta \pi_{A_3}) \text{ is positive semidefinite}.
\end{align*}

So, suppose that there exists an $\alpha > 0$ such that $\pi_{A_1} + \pi_{A_2} - \alpha \pi_{A_3}$ is positive semidefinite. It follows that there exists an $\gamma \geq \alpha$ such that $\pi_{A_1} + \pi_{A_2} - \gamma \pi_{A_3}$ is rank $1$ and positive semidefinite.  Since $\pi_{A_1} + \pi_{A_2} - \gamma \pi_{A_3}$ is a rank-1 positive semidefinite element of $\spn(\MES)$, we may apply Proposition \ref{MESspan} to state that there exists a $k>0$ and coisometry $B:\s Y \to \s X$ such that
\[
\pi_{A_1} + \pi_{A_2} - \gamma \pi_{A_3} = k \pi_B
\]
Thus, we conclude that
\begin{align*}
    \Phi(\pi_{A_1} + \pi_{A_2} - \alpha \pi_{A_3}) &= 
    \Phi(\pi_{A_1} + \pi_{A_2} - \gamma \pi_{A_3} + (\gamma - \alpha) \pi_{A_3})
    \\ & =
    \Phi(k \pi_B + (\gamma - \alpha) \pi_{A_3})
    \\ & = 
    k \Phi(\pi_B) + (\gamma - \alpha)\Phi(\pi_{A_3})
\end{align*}
Thus, the desired condition holds with $\beta = \alpha$. 

\proofend

\section{Pairwise Orthogonal Coisometries}

We will say that two coisometries $A_1,A_2: \s Y \to \s X$ are \textit{orthogonal} if $A_1A_2^* = 0$.

\begin{lem} \label{lem1}
For coisometries $A_1,A_2: \s Y \to \s X$, the following are equivalent:
\begin{enumerate}
    \item $A_1A_2^* = 0$.
    \item $A_2A_1^* = 0$.
    \item The map $\smat{A_1\\A_2}:\s Y \to \s X \oplus \s X$ is a coisometry.
    \item $\alpha A_1 + \beta A_2$ is a coisometry for $\alpha,\beta \in \C$ whenever $|\alpha|^2 + |\beta|^2 = 1$.
    \item The images of $A_1^*, A_2^*$ are orthogonal subspaces of $\s Y$.
\end{enumerate}

\end{lem}

\textit{Proof:} $1 \iff 2$: It suffices to note that $(A_1A_2^*)^* = A_2A_1^*$.

$1 \implies 3$: Suppose that $A_2 A_1^* = A_1A_2^* = 0$.  We then compute
\begin{equation}
    \mat{A_1\\A_2}\mat{A_1\\A_2}^* = \mat{A_1A_1^* & A_1A_2^*\\ A_2A_1^* & A_2A_2^*} = \mat{\id_{\s X} & 0\\0&\id_{\s X}} = \id_{\s X \oplus \s X}.
\end{equation}

$3 \implies 1$: Suppose that $\smat{A_1\\A_2}$ is a coisometry.  Then
\begin{equation}
    \id_{\s X \oplus \s X} =
    \mat{\id_{\s X} & 0\\0&\id_{\s X}} =
    \mat{A_1\\A_2}\mat{A_1\\A_2}^* = \mat{A_1A_1^* & A_1A_2^*\\ A_2A_1^* & A_2A_2^*},
\end{equation}
from which we may conclude that $A_jA_j^* = \id_{\s X}$ (that is, $A_1,A_2$ are coisometries) and $A_1 A_2^* = 0$ (that is, $A_1$ and $A_2$ are orthogonal).

$1 \iff 4$: We compute
\begin{align} 
    (\alpha A_1 + \beta A_2)(\alpha A_1 + \beta A_2)^* &= 
    |\alpha|^2 A_1A_1^* + |\beta|^2 A_2A_2^* 
    \\ & \qquad + \alpha \bar \beta A_1 A_2^* + \bar \alpha \beta A_2 A_1^*
    \\ & = (|\alpha|^2 + |\beta|^2) \id_{\s X} + (\alpha \bar \beta A_1 A_2^*) + (\alpha \bar \beta A_1 A_2^*)^*. \label{isomeq}
\end{align}
If $A_1A_2^* = 0$, then we have $(\alpha A_1 + \beta A_2)(\alpha A_1 + \beta A_2)^* = (|\alpha|^2 + |\beta|^2) \id_{\s X}$, which is to say that $\alpha A_1 + \beta A_2$ is a coisometry whenever $|\alpha|^2 + |\beta|^2 = 1$.  That is, $1 \implies 4$.

Conversely, if 4 holds, then it must be that $(\alpha \bar \beta A_1 A_2^*) + (\alpha \bar \beta A_1 A_2^*)^* = 0$ for all choices of $\alpha,\beta$ with $|\alpha|^2 + |\beta|^2 = 1$.  Let $\gamma = \alpha \bar \beta$, and choose self-adjoint $H,K$ such that $A_1A_2^* = H + iK$.  We then compute
\begin{equation}
    (\gamma(H + iK)) + (\gamma(H + iK))^* = 
    (\gamma + \bar \gamma)H + i(\gamma - \bar \gamma)K
\end{equation}
and, by 4, we know that the above expression is zero for all $\gamma = \alpha \bar \beta$ with $\alpha, \beta \in \C$ satisfying $|\alpha|^2 + |\beta|^2 = 1$.

Setting $\alpha = \beta = \frac{1}{\sqrt{2}}$, we find $\gamma = \frac 12$ and we see that we must have $H = 0$.  Setting $\alpha = \frac{i}{\sqrt{2}}$ and $\beta = \frac{1}{\sqrt{2}}$, we find $\gamma = \frac i2$ and we see that we must have $K = 0$.  Since $H = 0$ and $K = 0$, we conclude that $A_1A_2^* = 0$.

The proof that $5 \iff 1$ is straightforward and therefore omitted.

\proofend

%
%

\begin{lem} \label{lem2}
For any two orthogonal coisometries $A_1,A_2: \s Y \to \s X$: there exist orthogonal coisometries $B_1,B_2 : \s Y \to \s X$ with $B_i \in \zeta(A_i)$ such that we have either
\begin{align}
\Phi\left(\sum_{i,j = 1}^2 a_{ij}\vecop(A_i)\vecop(A_j)^*\right) &= 
\sum_{i,j = 1}^2 a_{ij}\vecop(B_i)\vecop(B_j)^*
\quad \text{or} \label{mesc1} \\
\Phi\left(\sum_{i,j = 1}^2 a_{ij}\vecop(A_i)\vecop(A_j)^*\right) &= 
\sum_{i,j = 1}^2 a_{ij}\vecop(B_j)\vecop(B_i)^*. \label{mesc2}
\end{align}
\end{lem}

\textit{Proof:} To begin, select coisometries $\tilde B_1 \in \zeta(A_1), \tilde B_2 \in \zeta(A_2)$. Because $\Phi$ is invertible, $\tilde B_1$ and $\tilde B_2$ must be linearly independent. Let $V_1$ and $V_2$ denote the subspaces
\begin{align*}
    V_1 &= \spn\{\vecop(A_i)\vecop(A_j)^*: \ i,j=1,2
\},\\
    V_2 &= \spn\{\vecop(\tilde B_i)\vecop(\tilde B_j)^*: \ i,j=1,2
\}.
\end{align*}
We observe that the elements of the form $\pi_{\alpha A_1 + \beta A_2}$ span $V_1$ (for instance: we see that $\{\pi_{A_1},\pi_{A_2},\pi_{(A_1+A_2)/2},\pi_{(A_1+iA_2)/2}\}$ forms a spanning set), and by Proposition \ref{3} we see that the image of these elements must lie in $V_2$.  Thus, $\Phi$ takes the subspace $V_1$ to $V_2$.
 
We define the map $G:\C^{2 \times 2} \to \C^{2 \times 2}$ so that
\begin{align*}
G\left(\mat{a_{11} & a_{12}\\ a_{21} & a_{22}}\right) &= \mat{b_{11} & b_{12}\\ b_{21} & b_{22}} \iff 
\\
\Phi\left(\sum_{i,j = 1}^2 a_{ij} \vecop(A_i)\vecop(A_j)^*\right)
&= \sum_{i,j = 1}^2 b_{ij} \vecop(\tilde B_i)\vecop(\tilde B_j)^* \iff
\\
\Phi\left(\smat{\vecop(A_1) & \vecop(A_2)}\mat{a_{11} & a_{12}\\ a_{21} & a_{22}}\smat{\vecop(A_1)^* \\ \vecop(A_2)^*}\right)
&= \smat{\vecop(\tilde B_1) & \vecop(\tilde B_2)}\mat{b_{11} & b_{12}\\ b_{21} & b_{22}}\smat{\vecop(\tilde B_1)^* \\ \vecop(\tilde B_2)^*}
\end{align*}
By our definition of $\tilde B_1,\tilde B_2,$ we have $G(E_{11}) = E_{11}$, $G(E_{22}) = E_{22}$. 

We note that
\begin{equation}\label{rk1}
\pi_{\alpha A_1 + \beta A_2} = 
\frac 1m  \smat{\vecop(A_1) & \vecop(A_2)}\left(
\mat{\alpha\\ \beta}\mat{\alpha \\ \beta}^* 
\right)\smat{\vecop(A_1)^* \\ \vecop(A_2)^*}
\end{equation}
so that because $\Phi$ preserves $\MES$, $G$ maps the rank-1 positive semidefinite matrices in $\C^{2 \times 2}$ to the rank-1 positive semidefinite matrices in $\C^{2 \times 2}$.  Thus, $G: \C^{2 \times 2} \to \C^{2 \times 2}$ is a positive map.  That is: if $A \in \C^{2 \times 2}$ is positive semidefinite, then $G(A)$ is also positive semidefinite.

By \cite{PosMaps}, $G$ must be decomposable.  That is, $\Phi$ must have the form $G = G_1 + G_2$, where $G_1$ and $\tau \circ G_2 $ are completely positive and $\tau$ denotes the transpose map $A \mapsto A^T$. 

Let $J(G)$ denote the Choi matrix of $G$.  That is, 
\[
J(G) =  \sum_{i,j = 1}^2 E_{ij} \otimes G (E_{ij})  
= (\id_{\C^{2 \times 2}} \otimes G)(\vecop(\id_{\C^2})\vecop(\id_{\C^2})^*)).
\]

So far, we may deduce that
\[
J(G) = J(G_1) + J(G_2) = 
\mat{G(E_{11}) & G(E_{12})\\ G(E_{21}) & G(E_{22})} =
\left[
\begin{array}{cc|cc}
1&0 & g_{11} & g_{12}\\
0&0 & g_{21} & g_{22}\\
\hline
\bar{g}_{11} & \bar{g}_{21} & 0&0\\
\bar{g}_{12} & \bar{g}_{22} & 0&1
\end{array}
\right]
\]
Let $P_i,Q_i,R_i$ be defined so that
\[
J(G_1) = \mat{P_1 & Q_1\\Q_1^* & R_1}, \quad 
J(G_2) = \mat{P_2 & Q_2\\Q_2^* & R_2}
\]
Because $G_1,\tau \circ G_2$ are completely positive, their Choi matrices $J(G_1)$ and $J(\tau \circ G_2 )$ must be positive semidefinite.
Thus, $P_1,P_2,R_1,R_2 \succeq 0$ and
\begin{equation}\label{2dcp}
\mat{P_1 & Q_1\\Q_1^* & R_1} \succeq 0, \quad 
\mat{P_2^T & Q_2^T\\\bar Q_2 & R_2^T} \succeq 0
\end{equation}
Because $P_1,P_2$ are positive definite matrices satisfying $P_1 + P_2 = E_{11}$, we must have $P_i = p_i E_{11}$ with $p_1 + p_2 = 1$ and $p_i \geq 0$.  Similarly, $R_i = r_i E_{22}$ with $r_1 + r_2 = 1$ and $r_i \geq 0$.

By \eqref{2dcp}, we must have $Q_1 = q_1 E_{12}$ and $Q_2 = q_2 E_{21}$, where $|q_i| \leq \sqrt{p_ir_i}$ for $i=1,2$.

Now, since $G$ maps rank-1 positive semidefinite matrices to rank-1 positive semidefinite matrices (as stated below \eqref{rk1}), it must be that the matrix
\[
G\left( \mat{\alpha\\ \beta}\mat{\alpha \\ \beta}^*\right) = 
\mat{
|\alpha|^2 & \alpha \bar \beta q_1 + \bar \alpha \beta \bar q_2\\
\bar \alpha \beta \bar q_1 + \alpha \bar \beta q_2 & |\beta|^2
}
\]
is positive semidefinite with rank $1$. Thus, its determinant must be zero for all $\alpha, \beta$.  That is, we must have
\[
|\alpha|^2 |\beta|^2 = |\alpha \bar \beta q_1 + \bar \alpha \beta \bar q_2|^2
\]
for all choices $\alpha,\beta \in \C$. Taking $\alpha = e^{i \theta}$ and $\beta = 1$, we see that this implies
\[
1 = |e^{i \theta} q_1 + e^{-i\theta} q_2|^2, \quad \theta \in \R.
\]
Thus, for all values of $\theta$, 
\[
|e^{i \theta} q_1 + e^{-i\theta} q_2|^2 = 
|q_1|^2 + |q_2|^2 + 2 \re[q_1 \bar q_2e^{i \theta}]
\]
is constant.  This only occurs if $q_1 \bar q_2 = 0$, which is to say that $q_1 = 0$ or $q_2 = 0$.  That is, we must have
\[
J(G) = \left[
\begin{array}{cc|cc}
1&0 & 0&e^{-i \alpha}\\
0&0 & 0&0\\
\hline
0&0 & 0&0\\
e^{i \alpha}&0 & 0&1
\end{array}
\right] \quad \text{or} \quad 
J(G) = \left[
\begin{array}{cc|cc}
1&0 & 0&0\\
0&0 & e^{i \alpha}&0\\
\hline
0&e^{-i \alpha} & 0&0\\
0&0 & 0&1
\end{array}
\right].
\]
%
%
%
%
If $J(G)$ has the first form, then we have
\begin{align*}
\Phi\left(\sum_{p,q = 1}^2 a_{pq}\vecop(A_p)\vecop(A_q)^*\right) &= 
\sum_{p,q = 1}^2 a_{pq}e^{i\alpha(p-q)}\vecop(\tilde B_p)\vecop(\tilde B_q)^*
\\ & = 
\sum_{p,q = 1}^2 a_{pq}\vecop(B_p)\vecop(B_q)^*,
\end{align*}
where we have taken $B_1 := \tilde B_1$ and $B_2:= e^{i \alpha} \tilde B_2$.  That is, \eqref{mesc1} applies, as desired.

If $J(G)$ has the second form, then we have
%
%
\begin{align*}
\Phi\left(\sum_{p,q = 1}^2 a_{pq}\vecop(A_p)\vecop(A_q)^*\right) &= 
\sum_{p,q = 1}^2 a_{pq}e^{i\alpha(q-p)}\vecop(\tilde B_q)\vecop(\tilde B_p)^*
\\ & = 
\sum_{p,q = 1}^2 a_{pq}\vecop(B_q)\vecop(B_p)^*,
\end{align*}
where again, we have taken $B_1 := \tilde B_1$ and $B_2:= e^{i \alpha} \tilde B_2$.  That is, \eqref{mesc2} applies, as desired.

If \eqref{mesc1} applies, then since $\Phi$ preserves $\MES$, we can state by Lemmas \ref{PT}, \ref{PTMES}, and \ref{lem1} that for all $\alpha,\beta$ with $|\alpha|^2 + |\beta|^2 = 1$, we have
\begin{align*}
\id_{\s X} &= 
\tr_{\s Y}(\Phi[(\vecop(\alpha A_1 + \beta A_2)
\vecop(\alpha A_1 + \beta A_2)^*)])
\\ &=
\tr_{\s Y}[|\alpha|^2 \vecop(B_1)\vecop(B_1)^*
+ |\beta|^2\vecop(B_2)\vecop(B_2)^*
\\ & \qquad 
+ \alpha \bar \beta \vecop(B_1)\vecop(B_2)^*
+ \bar \alpha \beta \vecop(B_2)\vecop(B_1)^*]
\\ & =
|\alpha|^2B_1B_1^* + |\beta|^2B_2B_2^* + 
\alpha \bar \beta B_1B_2^* + \bar \alpha \beta B_2B_1^*
\\ & = 
(|\alpha|^2 + |\beta|^2)\id_{\s X} + 2\re[\alpha \bar \beta B_1B_2^*]
= \id_{\s X} + 2\re[\alpha \bar \beta B_1B_2^*].
\end{align*}
Thus, we have $\re[\alpha \bar \beta B_1B_2^*] = 0$ for all $\alpha,\beta \in \C$ with $|\alpha|^2 + |\beta|^2 = 1$, so that $B_1B_2^* = 0$.  Thus, $B_1$ and $B_2$ are orthogonal coisometries.

Similarly: if \eqref{mesc2} applies, then for all $\alpha,\beta$ with $|\alpha|^2 + |\beta|^2 = 1$, we have
\begin{align*}
\id_{\s X} &= 
\tr_{\s Y}(\Phi[(\vecop(\alpha A_1 + \beta A_2)
\vecop(\alpha A_1 + \beta A_2)^*)])
\\ &=
\tr_{\s Y}[|\alpha|^2 \vecop(B_1)\vecop(B_1)^*
+ |\beta|^2\vecop(B_2)\vecop(B_2)^*
\\ & \qquad 
+ \alpha \bar \beta \vecop(B_2)\vecop(B_1)^*
+ \bar \alpha \beta \vecop(B_1)\vecop(B_2)^*]
\\ & = \id_{\s X} + 2\re[\alpha \bar \beta B_2B_1^*].
\end{align*}
We again deduce that $B_1$ and $B_2$ are orthogonal coisometries.  In either case, we have found that $B_1$ and $B_2$ are orthogonal coisometries, as desired.

\proofend

We now consider maps over $P(\C^2)$ (which would more commonly be notated as $P\C^1$ in topological contexts).  We will denote by $[\alpha,\beta]$ the equivalence class of $(\alpha,\beta) \in \C^2$.  That is, 
\[
[\alpha,\beta] = \{k (\alpha, \beta) : k \in \C\} \in P(\C^2)
\]

For coisometries $A,B$ satisfying $AB^* = 0$, we can select coisometries $C \in \zeta([A])$ and $D \in \zeta([B])$. For such an $A,B,C,D$, we then define $f : P(\C^2) \to P(\C^2)$ so that for $(\gamma,\delta) \in f([\alpha,\beta])$, we have
\[
\zeta([\alpha A + \beta B]) = [\gamma C + \delta D]
\]
That is, if $\iota_{A,B}$ denotes the homeomorphism $[\alpha,\beta] \in P(\C^2) \mapsto [\alpha A + \beta B] \in P(\spn\{A,B\})$, then 
\[
f_{A,B,C,D} = \iota_{C,D}^{-1} \circ \zeta \circ \iota_{A,B}
\]
so that $f$ is a well-defined, continuous map.

\begin{lem} \label{fcontlem}
For any coisometries $A,B$ satisfying $AB^* = 0$ and coisometries $C \in \zeta([A])$ and $D \in \zeta([B])$, the map $f:P(\C^2) \to P(\C^2)$ described above is continuous.
\end{lem}

\textit{Proof:} Let $\iota_{A,B}$ denote the homeomorphism $[\alpha,\beta] \in P(\C^2) \to [\alpha A + \beta B] \in P(\spn(A,B))$, and similarly define $\iota_{C,D}$.  We note that $\pi$ induces a homeomorphism between $P(\spn(A,B))$ and the image of $P(\spn(A,B))$ in $\MES_{\s X, \s Y}$, and a similar observation applies to $\spn(C,D)$.  Thus, we see that $f$ can be written as
\[
f = \iota_{C,D}^{-1} \circ \pi|_{P(\spn(C,D))}^{-1} \circ \Phi \circ \pi|_{P(\spn(A,B))} \circ \iota_{A,B}
\]
and is thus the composition of continuous functions. 

\proofend

Since $P(\C^2)$ is homeomorphic to $S^2$ (the two-dimensional sphere), we can consider the degree of the map $f_{A,B}:P(\C^2) \to P(\C^2)$.  Since $f$ is a homeomorphism, we must have $\deg(f) \in \{1,-1\}$.

\begin{lem} \label{deglem}

For a fixed $\MES$-preserving $\Phi$, there exists an $\epsilon \in \{1,-1\}$ such that for all choices of $A,B,C,D \in \coisom(\s Y, \s X)$ with $AB^* = 0$, $C\in \zeta([A])$, and $D \in \zeta([B])$, we have $\deg(f_{A,B}) = \epsilon$.

\end{lem}


\textit{Proof:} Let $f_0$ be the map $f$ corresponding to the coisometries $A_0,B_0,C_0,D_0$, and let $f_1$ be the map corresponding to the coisometries $A_1,B_1,C_1,D_1$.  In order to show that $f_0$ and $f_1$ have the same degree, it suffices to show that these maps are homotopic.

By Lemma \ref{lem1}, the set of pairs $(A,B)$ of coisometries such that $AB^* = 0$ is path connected.  Thus, there exist paths $A,B:[0,1] \to \coisom(\s Y, \s X)$ such that $A(0) = A_0$, $A(1) = A_1$, $B(0) = B_0$, $B(1) = B_1$, and for all $t \in [0,1]$ we have $A(t)B^*(t) = 0$.

Let $\Gamma(t) = \zeta( [A(t)])$ and $\Delta(t) = \zeta([B(t)])$.   $\Gamma, \Delta$ are paths in $P(\coisom(\s Y, \s X))$ connecting $\zeta([A_0])$ to $\zeta([A_1])$ and $\zeta([B_0])$ to $\zeta([B_1])$.

If we consider the tautological vector bundle over $P(\coisom(\s Y, \s X))$, it suffices to find a non-vanishing (constant magnitude) vector-field over the image $\Gamma$ connecting $(C_0, \zeta([A_0]))$ to $(C(1),\zeta([A_1]))$ and a non-vanishing vector-field over the image of $\Delta$ connecting $(D_0, \zeta([B_0]))$ to $(D_1,\zeta([B_1]))$.  These vector fields give us paths $C(t),D(t)$ in $\coisom(\s Y, \s X)$.

We can construct such a vector field as follows.  To begin, select a continuous $v,w:[0,1] \to \s X \otimes \s Y$ such that for all $t \in [0,1]$, $\pi_{\zeta([A_t])}v(t) \neq 0$ and $\pi_{\zeta([B_t])}w(t) \neq 0$.  With that, take $c(t) = \vecop^{-1}(\pi_{\zeta([A_t])}v(t))$.  Finally, define $C(t) = e^{i(p + qt)}\frac{c(t)}{\|c(t)\|}$, where $p,q \in \R$ are chosen so that $C(0) = C_0$ and $C(1) = C_1$.  
Similarly, take $d(t) = \vecop^{-1}(\pi_{\zeta([B_t])}w(t))$.  Define $D(t) = e^{i(r + st)}\frac{d(t)}{\|d(t)\|}$, where $r,s \in \R$ are chosen so that $D(0) = D_0$ and $D(1) = D_1$.

With that, we see that $f_t$ for $t \in [0,1]$ is a homotopy of the maps $f_0$ and $f_1$. 

\proofend

Combining Lemmas \ref{lem2} and \ref{deglem} lets us deduce the following.

\begin{prop} \label{semiprop}
Let $A_1,A_2$ be coisometries such that $A_1A_2^* = 0$.  Then there exist coisometries $B_1 \in \zeta([A_1]), B_2 \in \zeta([A_2])$ such that for any $\alpha,\beta \in \C$ with $|\alpha|^2 + |\beta|^2 = 1$, we have
\[
\zeta([\alpha A_1 + \beta A_2]) = [\alpha B_1 + \beta B_2]
\]
in the case that $\epsilon = +1$, or
\[
\zeta([\alpha A_1 + \beta A_2]) = [\overline{\alpha} B_1 + \overline{\beta}  B_2]
\]
in the case that $\epsilon = -1$.
\end{prop}

\textit{Proof:} In the case that \eqref{mesc1} holds for some coisometries $A_1,A_2,B_1,B_2$, we find that $\deg(f_{A_1,A_2,B_1,B_2}) = +1$, so that we have $\epsilon = +1$ globally. Applying \eqref{mesc1}, we find that
\begin{align*}
\Phi(\pi_{\alpha A_1 + \beta A_2}) &= 
\Phi\Big[\frac 1m(|\alpha|^2 \vecop(A_1)\vecop(A_1)^* + |\beta|^2 \vecop(A_2)\vecop(A_2)^*
\\ & \qquad +
\alpha \bar \beta \vecop(A_1)\vecop(A_2)^* + \bar \alpha \beta \vecop(A_2) \vecop(A_1)^*)
\Big]
\\ & = 
\frac 1m(|\alpha|^2 \vecop(B_1)\vecop(B_1)^* + |\beta|^2 \vecop(B_2)\vecop(B_2)^*
\\ & \qquad +
\alpha \bar \beta \vecop(B_1)\vecop(B_2)^* + \bar \alpha \beta \vecop(B_2) \vecop(B_1)^*)
\\ & = \pi_{\alpha B_1 + \beta B_2}
\end{align*}
so that we indeed have $\zeta([\alpha A_1 + \beta A_2]) = [\alpha B_1 + \beta B_2]$.

Similarly: in the case that \eqref{mesc2} holds for some coisometries $A_1,A_2,B_1,B_2$, we find that $\deg(f_{A_1,A_2,B_1,B_2}) = -1$, so that we have $\epsilon = -1$ globally. Applying \eqref{mesc2}, we find that
\begin{align*}
\Phi(\pi_{\alpha A_1 + \beta A_2}) &= 
\Phi\Big[\frac 1m(|\alpha|^2 \vecop(A_1)\vecop(A_1)^* + |\beta|^2 \vecop(A_2)\vecop(A_2)^*
\\ & \qquad +
\alpha \bar \beta \vecop(A_1)\vecop(A_2)^* + \bar \alpha \beta \vecop(A_2) \vecop(A_1)^*)
\Big]
\\ & = 
\frac 1m(|\alpha|^2 \vecop(B_1)\vecop(B_1)^* + |\beta|^2 \vecop(B_2)\vecop(B_2)^*
\\ & \qquad +
\alpha \bar \beta \vecop(B_2)\vecop(B_1)^* + \bar \alpha \beta \vecop(B_1) \vecop(B_2)^*)
\\ & = \pi_{\bar \alpha B_1 + \bar \beta B_2}
\end{align*}
so that we indeed have $\zeta([\alpha A_1 + \beta A_2]) = [\bar \alpha B_1 + \bar\beta B_2]$.

\proofend

\begin{rmk} 

Rather than using the degree of the maps $f_{A,B}$ to discriminate between the possible forms of $\zeta$, we could also use the $\det J(G)$, where $J(G)$ denotes the Choi matrix used in the proof of Lemma \ref{lem2}. In particular: in the $\epsilon = +1$ case we would compute $\det J(G) = 0$, and in the $\epsilon = -1$ case we would compute $\det J(G) = -1$. By the continuity of the determinant, one may argue that given an invertible $\MES$-preserving map $\Phi$, we must either have $\det J(G) = 0$ for all constructions of $G$, or $\det J(G) = -1$ for all constructions of $G$.

\end{rmk}

\section{Constructing an Extension}

To begin, we generalize Proposition \ref{semiprop}.

\begin{prop}\label{semipropgen}
Let $A_1,\dots,A_k: \s Y \to \s X$ be a collection of mutually orthogonal coisometries.  Then there exist mutually orthogonal coisometries $B_1,\dots,B_k$ such that for any $\alpha_1,\dots,\alpha_k \in \C$ with $|\alpha_1|^2 + \cdots + |\alpha_k|^2 = 1$, we have
\[
\zeta(\left[\alpha_1 A_1 + \cdots + \alpha_k A_k\right]) = 
\left[\alpha_1 B_1 + \cdots + \alpha_k B_k\right]
\]
in the case that $\epsilon = +1$, or 
\[
\zeta(\left[\alpha_1 A_1 + \cdots + \alpha_k A_k\right] )= 
\left[\bar \alpha_1 B_1 + \cdots + \bar \alpha_k B_k\right]
\]
in the case that $\epsilon = -1$.
\end{prop}

\textit{Proof:} In the $\epsilon = 1$ case, we can apply Lemma \ref{lem2} to show that there exist coisometries $B_1,\dots,B_k$ and $\theta_{pq} \in \R$ for $p,q = 1,\dots,k$ such that $\theta_{1,q} = 0$ for all $q$, $\theta_{qp} = -\theta_{pq}$, and we have
\[
\Phi\left( 
\sum_{p,q = 1}^k a_{pq} \vecop(A_p) \vecop(A_q)^*
\right) = 
\sum_{p,q = 1}^k a_{pq}e^{i\theta_{pq}} \vecop(B_p)\vecop(B_q)^*.
\]
Because $\Phi$ is $\MES$-preserving, the matrix $[e^{i\theta_{pq}}a_{pq}]_{p,q = 1}^k$ must be rank-one and positive semidefinite whenever $[a_{pq}]_{p,q=1}^n$ is rank-one and positive semidefinite.  By considering the case where $a_{pq} = 1$ for all $p,q$, 
we see that this can only occur if $M = [e^{i\theta_{ pq}}]_{p,q}$ has rank one.  If $M$ is a rank-one matrix whose first row and column are all $1$s, $M$ must be the matrix whose entries are all $1$s.  Thus, we have $e^{i\theta_{pq}} = 1$ for all $p,q$.  That is, we have
\[
\Phi\left( 
\sum_{p,q = 1}^k a_{pq} \vecop(A_p) \vecop(A_q)^*
\right) = 
\sum_{p,q = 1}^k a_{pq} \vecop(B_p)\vecop(B_q)^*,
\]
which is equivalent to the desired statement,
\[
\zeta(\left[\alpha_1 A_1 + \cdots + \alpha_k A_k\right]) = 
\left[\alpha_1 B_1 + \cdots + \alpha_k B_k\right].
\]

Similarly, in the $\epsilon = -1$ case, we can apply Lemma \ref{lem2} to show that there exist coisometries $B_1,\dots,B_k$ and $\theta_{pq} \in \R$ for $p,q = 1,\dots,k$ such that $\theta_{1,q} = 0$ for all $q$, $\theta_{qp} = -\theta_{pq}$, and we have
\[
\Phi\left( 
\sum_{p,q = 1}^k a_{pq} \vecop(A_p) \vecop(A_q)^*
\right) = 
\sum_{p,q = 1}^k a_{pq}e^{i\theta_{pq}} \vecop(B_q)\vecop(B_p)^*.
\]
Because $\Phi$ is $\MES$-preserving, the matrix $[e^{i\theta_{pq}}a_{pq}]_{p,q = 1}^k$ must be rank-one and positive semidefinite whenever $[a_{pq}]_{p,q=1}^n$ is rank-one and positive semidefinite. Applying the same analysis as above, we conclude $e^{i\theta_{pq}} = 1$ for all $p,q$.  That is, we have
\[
\Phi\left( 
\sum_{p,q = 1}^k a_{pq} \vecop(A_p) \vecop(A_q)^*
\right) = 
\sum_{p,q = 1}^k a_{pq} \vecop(B_q)\vecop(B_p)^*,
\]
which is equivalent to the desired statement,
\[
\zeta(\left[\alpha_1 A_1 + \cdots + \alpha_k A_k\right] )= 
\left[\bar \alpha_1 B_1 + \cdots + \bar \alpha_k B_k\right].
\]

\proofend

We will also need the following lemma, which is the polarization identity in one of its forms.

\begin{lem} \label{polarlem}
Given matrices $A,B \in \C^{p \times q}$, we have
\[
AB^* = \frac 14 \sum_{k=0}^3 i^k(A + i^k B) (A + i^k B)^*
\]
\end{lem}

\textit{Proof:}

We compute
\begin{align*}
    \sum_{k=0}^3 i^k (A + i^k B)(A + i^k B)^* &= \sum_{k=0}^3 i^k(AA^* + i^{-k} AB^* + i^k BA^* + BB^*)
    \\ & = 
    \sum_{k=0}^3 (i^k AA^* +  AB^* + i^{2k} BA^* + i^k BB^*)
    \\ & = 0 + 4AB^* + 0 + 0
\end{align*}
Thus, we have $AB^* = \frac 14 \sum_{k=0}^3 i^k(A + i^k B) (A + i^k B)^*$ as desired.

\proofend

\begin{prop} \label{ext}
Suppose that $\s Y = \s X^{k} := \s X \oplus \cdots \oplus \s X$.
Let $\Phi: \spn(\MES_{\s X,\s Y}) \to \spn(\MES_{\s X,\s Y})$ denote an invertible linear map such that $\Phi(\MES_{\s X, \s Y}) \subset \MES_{\s X, \s Y}$.  
Then $\Phi$ can be extended to a map $\tilde \Phi:\spn(\MES_{\s Y,\s Y}) \to \spn(\MES_{\s Y,\s Y})$ which satisfies $\tilde\Phi(\MES_{\s Y,\s Y}) \subset \MES_{\s Y,\s Y}$.
\end{prop}

\textit{Proof:} Let $\zeta:P(\coisom(\s Y, \s X)) \to P(\coisom(\s Y, \s X))$ be the map induced by $\Phi$, as defined in Section 4.
Let $U:\s Y \to \s Y$ be a unitary transformation given by $U = [U_j]_{j=1}^k$ (where $U_j:\s Y \to \s X$).  Note that $U_p,U_q$ are coisometries satisfying $U_pU_q^* = 0$ for $1 \leq p,q \leq k$. 

Define $\tilde \Phi:\MES_{\s Y, \s Y} \to \MES_{\s Y, \s Y}$ as follows: For $j = 1,\dots,k$ let $V_j \in \zeta([U_j])$ be such that $\zeta([\alpha U_p + \beta U_q]) = [\alpha V_p + \beta V_q]$ or $\zeta([\alpha U_p + \beta U_q]) = [\bar \alpha V_p + \bar \beta V_q]$ for $1 \leq p,q \leq k$, as guaranteed by Proposition \ref{semipropgen}.  Then, we take

\begin{align*}
    \tilde \Phi(&\pi_{[U_j]_j}) = 
    \tilde \Phi \left(
    \frac 1{km} [\vecop(U_p)\vecop(U_q)^*]_{p,q = 1}^k \right)
    \\ & = 
    \frac 1{km}[\vecop(V_p)\vecop(V_q)^*]_{p,q = 1}^k
\end{align*}
Equivalently, we have defined $\tilde \Phi$ so that
\[
\tilde \Phi(\pi_{[U_j]_j}) = \pi_{[V_j]_j}.
\]

%
%
By Proposition \ref{semipropgen}, $[V_j]_{j=1}^k$ is a coisometry, so that $\pi_{[V_j]_j} \in \MES_{\s Y, \s Y}$. Thus, we see from the above that $\tilde \Phi$ preserves $\MES_{\s Y, \s Y}$.  It remains to be shown, however, that $\tilde \Phi$ as defined above is a linear map.

In the case of $\epsilon = 1$, Lemma \ref{polarlem} yields
\begin{align*}
\vecop(U_p)\vecop(U_q)^* &= \frac 14\sum_{k=0}^3 i^k \vecop(U_p + i^k U_q)\vecop(U_p + i^k U_q)^*\\
\Phi[\vecop(U_p)\vecop(U_q)^*] &= \frac 14\sum_{k=0}^3 i^k \Phi[\vecop(U_p + i^k U_q)\vecop(U_p + i^k U_q)^*]\\
&= \frac 14\sum_{k=0}^3 i^k \vecop(V_p + i^k V_q)\vecop(V_p + i^k V_q)^*
\\ & = \vecop(V_p) \vecop(V_q)^*
\end{align*}
which means that our extension can be written as
\begin{align*}
    \tilde \Phi(\pi_{[U_j]_j}) &= 
    \tilde \Phi \left(
    \frac 1{km} [\vecop(U_p)\vecop(U_q)^*]_{p,q = 1}^k \right)
    \\ & =
    \frac 1{km} [\Phi(\vecop(U_p)\vecop(U_q)^*)]_{p,q = 1}^k.
\end{align*}
So, for every $M \in \MES_{\s Y, \s Y} \subset L(\s Y \otimes \s Y) = L(\C^k) \otimes L(\s X \otimes Y)$, we have
$\tilde \Phi(M) = (\id_{\C^{k \times k}}\otimes \Phi(M))$. That is, $\tilde \Phi|_{\MES_{\s Y, \s Y}}$ is the restriction of a linear map, and is therefore linear.

In the case of $\epsilon = -1$, Lemma \ref{polarlem} yields
\begin{align*}
\vecop(U_p)\vecop(U_q)^* &= \frac 14\sum_{k=0}^3 i^k \vecop(U_p + i^k U_q)\vecop(U_p + i^k U_q)^*\\
\Phi[\vecop(U_p)\vecop(U_q)^*] &= \frac 14\sum_{k=0}^3 i^k \Phi[\vecop(U_p + i^k U_q)\vecop(U_p + i^k U_q)^*]\\
&= \frac 14\sum_{k=0}^3 i^k \vecop(V_p + i^{-k}V_q)\vecop(V_p +  i^{-k}V_q)^*]\\
&= \frac 14\sum_{k=0}^3 i^k \left(i^{-k}\vecop(i^k V_p + V_q)\right)\left(i^{-k}\vecop(i^k V_p +  V_q)^*\right)]\\
&= \frac 14\sum_{k=0}^3 i^k \vecop(V_q + i^k V_p)\vecop(V_q + i^k V_p)^*]
\\ & = \vecop(V_q) \vecop(V_p)^*
\end{align*}
which means that our extension can be written as
\begin{align*}
    \tilde \Phi(\pi_{[U_j]_j}) &= 
    \tilde \Phi \left(
    \frac 1{km} [\vecop(U_p)\vecop(U_q)^*]_{p,q = 1}^k \right)
    \\ & =
    \frac 1{km} [\Phi(\vecop(U_q)\vecop(U_p)^*)]_{p,q = 1}^k
\end{align*}
So, for every $M \in \MES_{\s Y, \s Y} \subset L(\s Y \otimes \s Y) = L(\C^k) \otimes L(\s X \otimes \s Y)$, we have
$\tilde \Phi(M) = (\tau_{\C^{k \times k}}\otimes \Phi)(M)$. That is, $\tilde \Phi|_{\MES_{\s Y, \s Y}}$ is the restriction of a linear map, and is therefore linear.

\proofend

\section{The preservers of MES}

We now consider the possible forms of this extension, using Theorem 1 from \cite{Poon15}, which is to say Theorem \ref{Poonthm1} from the introduction. We note in particular that if $\Phi:\spn(\MES_{\s Y, \s Y}) \to \spn(\MES_{\s Y, \s Y})$ is invertible, then it must be of the form (1) or (2).

Let $\ad_U$ denote the map $\ad_U : A \mapsto UAU^{-1}$.
We make the following observation regarding these extended maps:

\begin{lem} \label{preadlem}

If $B \in L(\s X \otimes \s Y)$ satisfies $BM = MB$ for all $M \in \MES_{\s X, \s Y}$, then $B = c \id_{\s X \otimes \s Y}$ for some $c \in \C$.

\end{lem}

\textit{Proof:} For every maximally entangled vector $u \in \s X \otimes \s Y$, we have
\begin{align*}
Buu^* = uu^*B 
\implies (Buu^*)u = (uu^*B)u \implies 
Bu = (u^*Bu)u .
\end{align*}
That is, every maximally entangled vector $u$ is an eigenvector of $B$ with associated eigenvalue $u^*Bu$.
Thus, for every maximally entangled vector $u$, $u^*Bu$ is an eigenvalue of $B$.  Because the set of maximally entangled vectors is connected and the map $u \mapsto u^*Bu$ is continuous, the set
\[
\Omega = \{
u^*Bu : u \in \s X \otimes \s Y \text{ is maximally entangled}
\}
\]
is connected.  Also, $\Omega$ is a subset of the spectrum of $B$, which is a finite set.  Since $\Omega$ is connected and finite, it must be a singleton; thus there is a $c \in \C$ such that $\Omega = \{ c \}$.  Consequently, $Bu = c u$ for all maximally entangled vectors $u$.
Since the maximally entangled vectors span $\s X \otimes \s Y$, we may conclude that $B = c {\id_{\s X \otimes \s Y}}$ as desired. \hfill $\Box$

\begin{lem} \label{adlem}
If $U_1,U_2\in L(\s X)$ and $V_1,V_2 \in L(\s Y)$ are unitary operators such that 
\[
\ad_{U_1 \otimes V_1}(M) = \ad_{U_2 \otimes V_2}(M) \qquad \text{for all } M \in \MES_{\s X, \s Y}
\]
then we must have $U_1 \otimes V_1 = c\,U_2 \otimes V_2$ for some $c \in \C$. 
\end{lem}

%
%
%
%

\textit{Proof:} We may rewrite the above condition as
\[
(U_2 \otimes V_2)^*(U_1 \otimes V_1) M= M (U_2 \otimes V_2)^*(U_1 \otimes V_1) \qquad \text{for all } M \in \MES_{\s X, \s Y}.
\]
Lemma \ref{preadlem} now gives $(U_2 \otimes V_2)^*(U_1 \otimes V_1) = c\id_{\s X \otimes \s Y}$ for some $c \in \C$.  \proofend

\begin{rmk}
From here on, we suppose that $\s Y = \s X^k$, for some integer $k \geq 2$.  We also identify $\s Y = \s X^k$ with $\C^k \otimes \s X$, where it is convenient to do so.
\end{rmk}

We now observe that $\tilde \Phi$ has the following properties:

\begin{lem} \label{PQlem}
Define
\begin{align*} 
P_j &= [\id_{\C^k} - 2E_{jj}]\otimes \id_{\s X} \in L(\C^k \otimes \s X) = L(\s Y), 
\\
Q_{pq} &=  \left(T_{pq}\otimes \id_{\s X}\right)\otimes \id_{\s Y} \in L((\C^k \otimes \s X) \otimes \s Y) = L(\s Y \otimes \s Y),
\end{align*}
where $T_{pq}$ is the permutation matrix corresponding to the transposition of the $p,q$ entries of a vector in $\C^k$.

Then $\tilde \Phi \circ \ad_{ P_j \otimes \id_{\s Y}} = \ad_{P_j \otimes \id_{\s Y}} \circ \tilde \Phi$, and $\tilde \Phi \circ \ad_{Q_{pq}} = \ad_{Q_{pq}} \circ \tilde \Phi$.

\end{lem}

\textit{Proof:} In order to avoid cumbersome notation, we prove the result for $P_j$ with $j = 1$ and for $Q_{pq}$ with $p = 1, q = 2$. However, the same proof can be applied for arbitrary indices.

Let $M$ be an element of $\MES_{\s Y, \s Y}$, with
\[
M = [M_{p,q}]_{p,q=1}^k, \quad M_{pq} \in L(\s X \otimes \s Y).
\]
We compute
\begin{align*}
\ad_{P_1 \otimes \id_{\s Y}}(\tilde\Phi(M)) &= \tilde \Phi(\ad_{P_1 \otimes \id_{\s Y}}(M)) 
\\ & = 
\begin{cases}
\mat{\Phi(M_{11}) & -\Phi(M_{12}) & -\Phi(M_{13}) & \cdots \\ 
-\Phi(M_{21}) & \Phi(M_{22}) & \Phi(M_{23}) & \cdots \\
-\Phi(M_{31}) & \Phi(M_{32}) & \Phi(M_{33}) & \cdots \\
\vdots & \vdots & \vdots & \ddots} & \epsilon = 1\\
\qquad \\
\mat{\Phi(M_{11}) & -\Phi(M_{21}) & -\Phi(M_{31}) & \cdots \\ 
-\Phi(M_{12}) & \Phi(M_{22}) & \Phi(M_{32}) & \cdots \\
-\Phi(M_{13}) & \Phi(M_{23}) & \Phi(M_{33}) & \cdots \\
\vdots & \vdots & \vdots & \ddots} & \epsilon = -1
\end{cases}.
\end{align*}
And similarly,
\begin{align*}
\ad_{Q_{12}}(\tilde\Phi(M)) &= \tilde \Phi(\ad_{Q_{12}}(M)) 
\\ & = 
\begin{cases}
\mat{\Phi(M_{22}) & \Phi(M_{21}) & \Phi(M_{23}) & \cdots \\ 
\Phi(M_{12}) & \Phi(M_{11}) & \Phi(M_{13}) & \cdots \\
\Phi(M_{32}) & \Phi(M_{31}) & \Phi(M_{33}) & \cdots \\
\vdots & \vdots & \vdots & \ddots} & \epsilon = 1\\
\qquad \\
\mat{\Phi(M_{22}) & \Phi(M_{12}) & \Phi(M_{32}) & \cdots \\ 
\Phi(M_{21}) & \Phi(M_{11}) & \Phi(M_{31}) & \cdots \\
\Phi(M_{23}) & \Phi(M_{13}) & \Phi(M_{33}) & \cdots \\
\vdots & \vdots & \vdots & \ddots} & \epsilon = -1
\end{cases}.
\end{align*}

\proofend

\begin{lem} \label{commcomplem}

Fix any unitary $U,V \in L(\s Y)$, and let $\mb S: L(\s Y \otimes \s Y)\to L(\s Y \otimes \s Y)$ be the switch-operator, defined by $A \otimes B \mapsto B \otimes A$.  Then for unitary $A \in L(\s Y)$, we have
\begin{align*}
    \mb S(\pi_A) &= \pi_{A^T},\\
    (\pi_A)^T &= \pi_{\overline{A}},\\
    \ad_{U \otimes V}(\pi_A) &= \pi_{UAV^T}.
\end{align*}

\end{lem}

\begin{lem} \label{mesycomlem}
We have $\s Y = \s X^k$, with $k \geq 2$.
Suppose that $\Psi:\spn(\MES_{\s Y, \s Y})\to \spn(\MES_{\s Y, \s Y})$ is a linear, invertible, $\MES$-preserving map that commutes with $\ad_{P_j \otimes \id_{\s Y}}$ for all $j$.  Then we must have
\[
\Psi(M) = (U \otimes V)M^\sigma(U\otimes V)^*,
\]
where $\sigma$ is either the transpose or identity map. Moreover, either $U$ must be block diagonal (in other words, each copy of $\s X\subset \s Y$ is an invariant subspace of $U$), or $k=2$ and $U=\begin{bmatrix} 0 & U_{12} \cr U_{21} & 0 \end{bmatrix}$ with $U_{12} , U_{21} \in L(\s X)$ unitary.

\end{lem}

\textit{Proof:} Suppose that $\Psi$ is a map of the form (1) from Theorem \ref{Poonthm1}, with $\sigma$ the identity map.  That is, $\Psi = \ad_{U \otimes V}$.  Since $\Psi$ commutes with $\ad_{P_j \otimes \id_{\s Y}}$ for all $j$, we have $\ad_{{(P_j \otimes \id_{\s Y})}(U \otimes V)} = \ad_{(U \otimes V)(P_j \otimes \id_{\s Y})}$. By Lemma \ref{adlem}, we have $(P_jU) \otimes V =c (UP_j) \otimes V$ for some $c \in \C$.  That is, if $U = [U_{p,q}]_{p,q = 1}^k$, then we have in the case of $j = 2$
\[
\mat{U_{11} & U_{12} & U_{13} & \cdots \\ 
-U_{21} & -U_{22} & -U_{23} & \cdots \\
U_{31} & U_{32} & U_{33} & \cdots \\
\vdots & \vdots & \vdots & \ddots } \otimes V = 
c\, \mat{U_{11} & -U_{12} & U_{13} & \cdots \\ 
U_{21} & -U_{22} & U_{23} & \cdots \\
U_{31} & -U_{32} & U_{33} & \cdots \\
\vdots & \vdots & \vdots & \ddots } \otimes V.
\]
When $k>2$, we conclude that $U_{pq} = U_{qp} = 0$ whenever $p \neq q$, as desired. When $k=2$ we either have that $c=1$ and $U_{12}=U_{21}=0$, or that $c=-1$ and $U_{11} = U_{22} = 0$.

Similarly, suppose that $\Psi$ is a map of the form (1) with $\sigma$ the transpose map, which is to say that $\Psi = \ad_{U \otimes V} \circ \tau$. Since $\Psi$ commutes with $\ad_{P_j \otimes \id_{\s Y}}$, we have 
\begin{align*}
\ad_{P_j \otimes \id_{\s Y}} \circ \ad_{U \otimes V} \circ \tau &= 
\ad_{U \otimes V} \circ \tau \circ \ad_{P_j \otimes \id_{\s Y}} \implies\\
\ad_{P_j \otimes \id_{\s Y}} \circ \ad_{U \otimes V} \circ \tau &= 
\ad_{U \otimes V} \circ \ad_{P_j \otimes \id_{\s Y}} \circ \tau \implies\\
\ad_{P_j \otimes \id_{\s Y}} \circ \ad_{U \otimes V} &=
\ad_{U \otimes V} \circ \ad_{P_j \otimes \id_{\s Y}}.
\end{align*}
By Lemma \ref{adlem}, we have $(P_jU) \otimes V = c\,(UP_j) \otimes V$ for some $c \in \C$, which again allows us to conclude that $U$ is block-diagonal, or $k=2$ and $U=\begin{bmatrix} 0 & U_{12} \cr U_{21} & 0 \end{bmatrix}$.

Now, suppose that $\Psi$ is of the form (2) from Theorem \ref{Poonthm1}, with $\sigma$ the identity map. That is, suppose that $\Psi = \ad_{U \otimes V} \circ \mb S$; we wish to show that $\Psi$ must fail to commute with $\ad_{P_j \otimes \id_{\s Y}}$ for some $j$.  By Lemma \ref{commcomplem}, we find that for unitary $A \in L(\s Y)$, we have
\begin{align*}
    (\ad_{P_j \otimes \id_{\s Y}} \circ \Psi)(\pi_A) &= \pi_{P_jUA^TV^T}\\
    (\Psi \circ \ad_{P_j \otimes \id_{\s Y}})(\pi_A) &= \pi_{UA^TP_jV^T}.
\end{align*}
Thus, it suffices to find unitary $A \in L(\s Y)$ such that for some $j$, 
$P_jUA^T  \neq c\, UA^TP_j$ for all $c \in \C$. Thus, we see that one $M \in \MES_{\s Y, \s Y}$ for which $(\ad_{P_j \otimes \id_{\s Y}} \circ \Psi)(M) \neq (\Psi \circ \ad_{P_j \otimes \id_{\s Y}})(M)$ is given by $M = \pi_A$, where $A$ is chosen so that $P_jUA^T$ and $UA^TP_j$ are not multiples.  For instance, for $j=1$, we can take
\[
UA^T = \frac 1{\sqrt{2}}\mat{\id_{\s X} & \id_{\s X}\\ \id_{\s X} & -\id_{\s X}
} \oplus \id_{\s X^{k-2}}.
\]

Similarly, suppose that $\Psi$ is of the form (2) from Theorem \ref{Poonthm1}, with $\sigma$ the transpose map. That is, suppose that $\Psi = \ad_{U \otimes V} \circ \tau \circ \mb S$.  By Lemma \ref{commcomplem}, we find that for unitary $A \in L(\s Y)$, we have
\begin{align*}
    (\ad_{P_j \otimes \id_{\s Y}} \circ \Psi)(\pi_A) &= \pi_{P_jUA^*V^T}\\
    (\Psi \circ \ad_{P_j \otimes \id_{\s Y}})(\pi_A) &= \pi_{UA^*P_jV^T}.
\end{align*}
Thus, it suffices to find unitary $A \in L(\s Y)$ such that for some $j$,
$    P_jUA^*  \neq c\, UA^*P_j $ for all $c \in \C $.
Thus, we see that one $M \in \MES_{\s Y, \s Y}$ for which $(\ad_{P_j \otimes \id_{\s Y}} \circ \Psi)(M) \neq (\Psi \circ \ad_{P_j \otimes \id_{\s Y}})(M)$ is given by $M = \pi_A$, where $A$ is chosen so that $P_jUA^*$ and $UA^*P_j$ are not multiples.  For instance, we can take
\[
UA^* = \frac 1{\sqrt{2}}\mat{\id_{\s X} & \id_{\s X}\\ \id_{\s X} & -\id_{\s X}
} \oplus \id_{\s X^{k-2}}.
\]

\proofend

\begin{lem} \label{mesycomlem2}
Let $U, V \in L(\s Y )$ be unitary with $U \in L(\s Y) = L(\C^k \otimes \s X)$ block diagonal. If $\Psi$ defined by
\[
\Psi(M) = (U \otimes V)M^\sigma(U\otimes V)^*, \ M \in \MES_{\s Y, \s Y},
\]
commutes with $\ad_{Q_{pq}}$ for $1 \leq p,q \leq k$, then $U$ either has the form
\[
U = \id_{\C^k} \otimes W,
\]
for some unitary $W:\s X \to \s X$, or $k=2$ and $U=\begin{bmatrix} 1 & 0 \cr 0 & -1 \end{bmatrix} \otimes W$ for some unitary $W:\s X \to \s X$.
\end{lem}

\textit{Proof:} As in the last proof (invoking Lemma \ref{preadlem}), we see that if $\Psi$ commutes with $\ad_{Q_{pq}}$, then we must have $Q_{pq}(U \otimes V) = c_{pq}\,(U \otimes V)Q_{pq}$ for some $c_{pq} \in \C$.  Let $\diag[U_j]_{j=1}^k$ denote the block-diagonal matrix with $U_j$ on the $j$th diagonal.  Then the above tells us that $U_p = c_{pq} U_q$, $U_q= c_{pq} U_p$ and $U_j=c_{pq} U_j$ for $j\neq p,q$. 

When $k>2$, this implies that $c_{pq}=1$ and $U_p=U_q$. As this holds for all pairs $p,q$, 
we may conclude that $U_j = U_{\ell}$ for all $1 \leq j,\ell \leq k$.  

When $k=2$, we get from $U_1= c_{12} U_2 = c_{12}^2 U_1$, that $c_{12} = \pm 1$. Thus either $U_1=U_2$ or $U_1=-U_2$.

\proofend

%
%
%

We can now finally prove our main result.

\textit{Proof of Theorem \ref{pres}:} By Proposition \ref{ext}, $\Phi$ has an extension $\tilde \Phi: L(\s Y \otimes \s Y) \to L(\s Y \otimes \s Y)$, as stated there.  By Lemma \ref{PQlem}, this extension commutes with $\ad_{P_j \otimes \id_{\s Y}}$ and $\ad_{Q_{pq}}$ for all $1 \leq j,p,q \leq k$.  By Lemma \ref{mesycomlem}, it follows that $\tilde \Phi$ is a map of the form $\tilde\Phi(M) = (\tilde U \otimes V)M^\sigma(\tilde U\otimes V)^*$, where $\sigma$ is either the transpose or identity map and with $\tilde U$ either block diagonal or, when $k=2$, of the form $\tilde U = \begin{bmatrix} 0 & U_{12} \cr U_{21} & 0 \end{bmatrix}$. 

We first show that we can discard the last possibility. Indeed, from the proof of Proposition \ref{ext} we obtain that
$$ \tilde\Phi (\begin{bmatrix} M_{11} & M_{12} \cr M_{21} & M_{22} \end{bmatrix} ) = \begin{bmatrix} \Phi(M_{11}) & \Phi(M_{12}) \cr \Phi (M_{21}) & \Phi (M_{22}) \end{bmatrix} \ {\text{or}}  \ \tilde\Phi (\begin{bmatrix} M_{11} & M_{12} \cr M_{21} & M_{22} \end{bmatrix} ) = \begin{bmatrix} \Phi(M_{11}) & \Phi(M_{21}) \cr \Phi (M_{12}) & \Phi (M_{22}) \end{bmatrix}. $$
If $\tilde U = \begin{bmatrix} 0 & U_{12} \cr U_{21} & 0 \end{bmatrix}$, then we obtain that for all $M=[M_{pq}]_{p,q=1}^2 \in \MES_{\s Y,\s Y}$ either
$$ \Phi (M_{11}) = (U_{12} \otimes V) M_{22} (U_{12} \otimes V)^*  \ {\text{or}}  \ \Phi (M_{11}) = (U_{12} \otimes V) M_{22}^T (U_{12} \otimes V)^*.$$ 
When $\dim X >1$ (and thus $\dim Y >2$) this leads to a contradiction as one can vary $M_{22}$ without changing $M_{11}$. 

If $k = 2$ and $\dim \s X = 1$, we have $\s X \otimes \s Y \cong \s Y = \C^2$, and $\Phi: L(\s X \otimes Y) \to L(\s X \otimes \s Y)$ defines a map on $2 \times 2$ matrices that preserves all rank-one orthogonal projections. This gives that $\Phi (X) = V X^\sigma V^*$, where $\sigma$ is either the transpose or the identity map and $V$ is unitary. Thus $\Phi$ is of the desired form.

Going back to the general case of $k\ge 2$, we can conclude that $\tilde U$ is block diagonal. Next, by Lemma \ref{mesycomlem2}, we find that 
\[
\tilde U = \id_{\C^k} \otimes U, \ \text{or} \ k=2 \ \text{and}\ \tilde U=\begin{bmatrix} 1 & 0 \cr 0 & -1\end{bmatrix} \otimes U ,
\]
for some unitary $U: \s X \to \s X$. We now compute for $M=[M_{pq}]_{p,q=1}^k \in \MES_{\s Y,\s Y}$ that 
\begin{align*}
\Phi(M_{11} )& = [e_1^T \otimes \id_{\s X} \otimes \id_{\s Y}]
\tilde \Phi(M) [e_1^T \otimes \id_{\s X} \otimes \id_{\s Y}]^*
\\ & = 
[e_1^T \otimes \id_{\s X} \otimes \id_{\s Y}]
(\tilde U \otimes V)M^\sigma(\tilde U\otimes V)^* [e_1 \otimes \id_{\s X} \otimes \id_{\s Y}]
\\ &= 
(U \otimes V)M_{11}^\sigma (U \otimes V)^*,
\end{align*}
as desired. 

\proofend

\end{document}